\documentclass{amsart}

\usepackage{amsmath,amssymb}
\usepackage{amsthm}
\usepackage[alphabetic]{amsrefs}
\usepackage{indentfirst}
\usepackage{mathrsfs}
\usepackage{graphicx}
\usepackage{hyperref}
\usepackage{xcolor}
\usepackage[margin=1.5in]{geometry}

\usepackage{todonotes}

\theoremstyle{plain}
\newtheorem{theorem}{Theorem}

\newtheorem{corollary}[theorem]{Corollary}

\theoremstyle{definition}

\theoremstyle{remark}

\DeclareMathOperator{\spanop}{span}

\IfFileExists{mathabx.sty}%
  {\DeclareFontFamily{U}{mathx}{\hyphenchar\font45}%
   \DeclareFontShape{U}{mathx}{m}{n}{<->mathx10}{}%
   \DeclareSymbolFont{mathx}{U}{mathx}{m}{n}%
   \DeclareFontSubstitution{U}{mathx}{m}{n}%
   \DeclareMathAccent{\widebar}{0}{mathx}{"73}%
}{%
  \PackageWarning{mathabx}{%
    Package mathabx not available, therefore\MessageBreak substituting
    widebar with overline\MessageBreak }%
  \newcommand{\widebar}[1]{\overline{#1}}%
}
\newcommand{\wb}[1]{\widebar{#1}}

\newcommand{\norm}[1]{\lVert#1\rVert}
\newcommand{\Norm}[1]{\left\lVert#1\right\rVert}
\newcommand{\average}[1]{\langle#1\rangle}

\title{On pointwise products of elliptic eigenfunctions}


\author[]{Jianfeng Lu}
\address[J. L.]{Department of Mathematics, Department of Physics, and Department of Chemistry,
Duke University, Box 90320, Durham NC 27708, USA}
\email{jianfeng@math.duke.edu}

\author[]{Stefan Steinerberger}
\address[S.S.]{Department of Mathematics, Yale University, New Haven, CT 06511, USA}
\email{stefan.steinerberger@yale.edu}

\thanks{J.L. is supported in part by the NSF National Science (DMS-1454939). S.S. is supported in part by the NSF (DMS-1763179) and the Alfred P. Sloan Foundation.}

\begin{document}

\begin{abstract} We consider eigenfunctions of Schr\"odinger
  operators 
  on a $d-$dimensional bounded domain $\Omega$ (or a $d-$dimensional
  compact manifold $\Omega$) with Dirichlet conditions. These
  operators give rise to a sequence of eigenfunctions
  $(\phi_n)_{n \in \mathbb{N}}$. We study the subspace of all
  pointwise products
$$ A_n = \mbox{span} \left\{ \phi_i(x) \phi_j(x): 1 \leq i,j \leq n\right\} \subseteq L^2(\Omega).$$
Clearly, that vector space has dimension $\mbox{dim}(A_n) = n(n+1)/2$. We prove that products $\phi_i \phi_j$ of eigenfunctions are simple in a certain sense: for any $\varepsilon > 0$, there exists a low-dimensional vector space $B_n$ that almost contains all products. More precisely, denoting the orthogonal projection $\Pi_{B_n}:L^2(\Omega) \rightarrow B_n$, we have
$$ \forall~1 \leq i,j \leq n~  \qquad \|\phi_i\phi_j - \Pi_{B_n}( \phi_i \phi_j) \|_{L^2} \leq \varepsilon$$
and the size of the space $\mbox{dim}(B_n)$ is relatively small
$$ \mbox{dim}(B_n) \lesssim    \Bigl(  \frac{1}{\varepsilon} \max_{1 \leq i \leq n} \norm{\phi_i}_{L^{\infty}} \Bigr)^d n.$$
In the generic delocalized setting, this bound grows linearly up to
logarithmic factors: pointwise products of eigenfunctions are
low-rank. This has implications, among other things, for the validity
of fast algorithms in electronic structure computations.
\end{abstract}

\maketitle

\section{Introduction}
Throughout this paper, we let $L$ be a uniformly elliptic second-order partial differential operator of Schr\"odinger type
$$ Lu = - \mbox{div}(a(x) \cdot \nabla u) + V(x)u$$
on a $d-$dimensional bounded domain $\Omega$ (or a smooth
$d-$dimensional compact manifold $\Omega$) with Dirichlet conditions
and smooth coefficients $a, V \in C^{\infty}$. This operator $L$ gives
rise to a sequence of eigenfunctions $L \phi_k = \lambda_k \phi_k$
which we assume to be $L^2-$normalized and which span
$L^2(\Omega)$. We ask a very simple question.
\begin{quote}
\textbf{Question.} What can be said about the function $\phi_{i}(x) \phi_{j}(x)$? 
\end{quote} 
Clearly, by $L^2-$orthogonality, the function
$\phi_{i}(x) \phi_{j}(x)$ has mean value 0 if $i\neq j$ but what else
can be said about its spectral resolution, for example, the size of
$\left\langle \phi_i \phi_j, \phi_k\right\rangle$?  There are very few
results overall; some results have been obtained in the presence of
additional structure assumptions on $\Omega$ connected to number
theory (see Bernstein \& Reznikoff \cite{bern}, Kr\"otz \& Stanton
\cite{krotz} and Sarnak \cite{Sarnak}). Already the simpler question
of understanding $L^2-$size of the product is highly nontrivial: a
seminal result of Burq-G\'{e}rard-Tzetkov \cite{burq} states
$$ \|\phi_{\mu}\phi_{\lambda} \|_{L^2} \lesssim \min( \sqrt{\lambda}, \sqrt{\mu})^{} \| \phi_{\lambda}\|_{L^2} \|\phi_{\mu}\|_{L^2}$$
on compact two-dimensional manifolds without boundary (this has been extended to higher dimensions \cite{blair, burq2}). A recent result of the second author \cite{stein}  (see also \cite{alex})
shows that one would generically (i.e., on typical manifolds in the presence of quantum chaos) expect $\phi_i(x)\phi_j(x)$ to be mainly supported at eigenfunctions having their eigenvalue close to $\max\left\{\lambda_i, \lambda_j\right\}$ and that deviation
from this phenomenon, as in the case of Fourier series on $\mathbb{T}$ for example, requires eigenfunctions to be strongly correlated at the wavelength in a precise sense.

In this paper, we ask the question on the numerical
rank of the space spanned by the pointwise products of eigenfunctions
$$ A_n = \mbox{span} \left\{ \phi_i(x) \phi_j(x): 1 \leq i,j \leq n\right\}.$$
This is a natural quantity for measuring the complexity of the products but also motivated by the density fitting approximation to the electron repulsion integral in the quantum chemistry literature. 
Given a set of eigenfunctions, the four-center two-electron repulsion integral
$$ (ij|kl) = \int_{\Omega \times \Omega}{\frac{\phi_i(x) \phi_j(x) \phi_k(y) \phi_l(y)}{|x-y|} dx dy}$$
is a central quantity in electronic structure theories. If we are
working with the first $n$ eigenfunctions
$(\phi_i)_{1 \leq i \leq n}$, then one has to evaluate $\mathcal{O}(n^4)$ integrals.

It has been empirically observed in the literature (see e.g.,
\cite{jian1} by the first author and Lexing Ying) that the space $A_n$
can in practice be very well approximated by another vector space
$B_n$ with $\mbox{dim}(B_n) \sim c \cdot n$, often referred as density
fitting.  This then drastically reduces the number of integrals in need of evaluation to
$\mathcal{O}(n^2)$. Our result is inspired by the empirical
success of density fitting and gives a mathematical justification.

\section{The results}
We state our main result: pointwise products of delocalized eigenfunctions are low-rank.
\begin{theorem}[Main Theorem] For any $\varepsilon>0, n \in \mathbb{N}$, there exists a space $B_n$ such that
$$ \forall~i,j \leq n~  \qquad \|\phi_i\phi_j - \Pi_{B_n}( \phi_i \phi_j) \|_{L^2} \leq \varepsilon$$
and such that $\emph{dim}(B_n)$ is relatively small
$$ \emph{dim}(B_n) \lesssim_{\Omega, a, V}    \Bigl(  \frac{1}{\varepsilon} \max_{1 \leq i \leq n} \norm{\phi_i}_{L^{\infty}} \Bigr)^d n.$$
\end{theorem}
We emphasize that for generic differential operators on generic domains, we expect the eigenfunctions to be delocalized and their $L^{\infty}-$norm to be governed by the random wave heuristic: in particular, we expect the maximum to grow logarithmically. Theorem 1 is thus optimal in the quantum chaotic setting (up to at most logarithmic factors). It is an interesting question whether the logarithm can be removed.
As a corollary of Theorem 1, we obtain unconditional nontrivial results in low dimensions $d=1,2$. 

\begin{corollary} We have
$$ \emph{dim}(B_n) \lesssim_{\Omega, a, V}  \begin{cases} \varepsilon^{-1} \cdot n \qquad &\mbox{if}~d=1 \\ \varepsilon^{-2} \cdot n^{3/2} \qquad &\mbox{if}~d=2. \end{cases}$$
\end{corollary}
The case $d=1$ is clearly optimal and well in line with classical intuition (the WKB asymptotic suggests that
things behave asymptotically like trigonometric functions which one would expect to behave like Fourier series whose products are low-rank).

\medskip 

In light of the  four-center two-electron repulsion integral, which is a somewhat better behaved quantity due to the smoothing effects of the potential, it is natural to look for a similar result in a function spaces
that captures the smoothing of the potential: since $|x-y|^{-1}$ is the fundamental solution of the Laplacian in $\mathbb{R}^3$, the appropriate space for physically relevant problems in $\Omega \subset \mathbb{R}^3$ is the homogeneous Sobolev space $H^{-1}$ norm (also known as the Coulomb norm or the Beppo-Levi norm) equipped with 
\begin{equation*}
  \norm{f}_{\dot{H}^{-1}}^2 = \bigl\langle  f - \wb{f}, \Delta^{-1} ( f - \wb{f}) \bigr\rangle, 
\end{equation*}
where $\wb{f}$ is the average of $f$, so that the above definition makes sense also for closed manifold. 

\begin{theorem}[Variant in $\dot H^{-1}$] For any $\varepsilon>0, n \in \mathbb{N}$, there exists a space $B_n$ such that
$$ \forall~i,j \leq n~  \qquad \|\phi_i\phi_j - \Pi_{B_n}( \phi_i \phi_j) \|_{\dot H^{-1}} \leq \varepsilon$$
such that $\emph{dim}(B_n)$ is relatively small
$$ \emph{dim}(B_n) \lesssim_{\Omega, a, V} 
    \Bigl( \frac{1}{\varepsilon} \max_{1\leq i \leq n} \norm{\phi_i}_{L^{\infty}} \Bigr)^{\frac{d}{2}}
    \sqrt{n}.$$
\end{theorem}

Note that this improves the bound in the $L^2$ case by essentially a square root. 
Here, we can obtain an unconditional result that improves on the
trivial bound for all dimensions $d \leq 6$. We also note that it
matches the linear scaling for $d=3$.

\begin{corollary}[Variant in $\dot H^{-1}$] We have
$$ \emph{dim}(B_n) \lesssim_{\Omega, a, V}  \varepsilon^{-\frac{d}{2}} n^{\frac{d+1}{4}}.$$
\end{corollary}

It is clear that the underlying low-rankness phenomenon is of great
relevance in a variety of settings (e.g., fast algorithms for
electronic structure theories \cite{jian1, jian2}, etc.); it is also
clear that similar arguments to the ones contained in this paper can
be used to obtain low-rankness results of a similar flavor in a
variety of function spaces. One of the most interesting questions is
whether it is possible to have $\mbox{dim}(B_n) \sim n$. Corollary 2
shows that this is the case for ordinary differential equations in
$L^2$; it would be interesting to have concrete examples and
counterexamples in higher dimensions -- unfortunately many of the
cases that can be made explicit (the torus $\mathbb{T}^d$, the sphere
$\mathbb{S}^d$) seem exceedingly non-generic.

\section{Proofs}

\begin{proof}[Proof of Theorem 1]
We first prove an elementary inequality: let $1 \leq i, j \leq n$ and
 $$-\mbox{div}(a(x)\nabla \phi_i) + V \phi_i = \lambda_i \phi_i  \qquad \mbox{and} \qquad  -\mbox{div}(a(x)\nabla \phi_j) + V \phi_j = \lambda_j \phi_j  $$ with the normalization $\norm{\phi_i}_{L^2} = \norm{\phi_j}_{L^2} = 1$, then we have 
  \begin{equation}\label{eq:quadest}
   \left\langle (-\mbox{div}(a(x)\nabla (\phi_i \phi_j)) + V \phi_i \phi_j, \phi_i \phi_j \right\rangle \lesssim_{a,V} (\lambda_n + 1) \max_{1 \leq i \leq n}{ \|\phi_i\|^2_{L^{\infty}}}.
 \end{equation}

 We separate the relevant quantity into two terms. The first one is fairly simple, we observe that $V \in C^{\infty}(\Omega)$ and thus, by $L^2-$normalization of the eigenfunctions,
\begin{align*}
   \left\langle V \phi_i \phi_j, \phi_i \phi_j \right\rangle &= \int_{\Omega}{V \phi_i^2 \phi_j^2 dx} \\
&\leq \|V\|_{L^{\infty}} \max_{1 \leq i \leq n}{ \|\phi_i\|^2_{L^{\infty}}} \lesssim_{V}  \max_{1 \leq i \leq n}{ \|\phi_i\|^2_{L^{\infty}}} 
\end{align*}
The second term can be rewritten using integration by parts as
\begin{align*}
  \left\langle (-\mbox{div}(a(x)\nabla ( \phi_i \phi_j)), \phi_i \phi_j  \right\rangle &= \int_{\Omega}{a(x) ( \nabla (\phi_i \phi_j) \cdot  \nabla (\phi_i \phi_j) ) dx}\\
&\leq \left( \max_{x \in \Omega}{a(x)} \right) \|  \nabla (\phi_i \phi_j)\|_{L^2}^2.
\end{align*}
The product rule and the triangle inequality yield
  \begin{equation*}
    \begin{aligned}
      \norm{ \nabla (\phi_i \phi_j) }_{L^2} & = \norm{ \phi_i \nabla \phi_j + \phi_j \nabla \phi_i}_{L^2} \\
      & \leq \norm{\nabla \phi_j}_{L^2} \norm{\phi_i}_{L^{\infty}} + \norm{\nabla \phi_i}_{L^2} \norm{\phi_j}_{L^{\infty}}.
      \end{aligned}
  \end{equation*}
However, we also have, with integration by parts,
\begin{align*}
 \norm{\nabla \phi_i}_{L^2}^2  &= \int_{\Omega}{ \nabla \phi_i \cdot \nabla \phi_i dx} \\
&\leq  \left(\min_{x \in \Omega}{a(x)}^{} \right)^{-1} \int_{\Omega}{ a(x) \nabla \phi_i \cdot \nabla \phi_i dx} \\
&= \left(\min_{x \in \Omega}{a(x)}^{} \right)^{-1} \int_{\Omega}{ - \mbox{div}(a(x) \nabla \phi_i)  \phi_i dx} \\
&=\left(\min_{x \in \Omega}{a(x)}^{} \right)^{-1} \int_{\Omega}{( \lambda_i \phi_i - V \phi_i)  \phi_i dx} \\
&\lesssim_{a} \lambda_n + \|V\|_{L^{\infty}}.
\end{align*}
This concludes the proof of \eqref{eq:quadest}.

\smallskip 

Expanding the product
\begin{equation*}
  \phi_i \phi_j = \sum_{k=1}^{\infty} \average{\phi_i \phi_j, \phi_k} \phi_k,  
\end{equation*}
we obtain, since $\phi_k$ are the eigenfunctions, 
\begin{align*}
   \bigl\langle (-\mbox{div}(a(x)\nabla (\phi_i \phi_j)) + V \phi_i \phi_j, \phi_i \phi_j \bigr\rangle &= \left\langle  \sum_{k=1}^{\infty} \lambda_k \average{\phi_i \phi_j, \phi_k} \phi_k,  \sum_{k=1}^{\infty} \average{\phi_i \phi_j, \phi_k} \phi_k\right\rangle\\
&=  \sum_{k=1}^{\infty} \lambda_k |\average{\phi_i \phi_j, \phi_k}|^2 .
\end{align*}
This shows, for any $r \in \mathbb{N}$,
  \begin{align*}
   \left\langle (-\mbox{div}(a(x)\nabla (\phi_i \phi_j)) + V \phi_i \phi_j, \phi_i \phi_j \right\rangle \geq  \lambda_r \Norm{\Pi_{> r}(\phi_i \phi_j) }_{L^2}^2,
  \end{align*}
  where we have introduced the short hand notation
  $\Pi_{> r} = \Pi_{\spanop\{\phi_k \mid k > r\}}$. Combined this with
  \eqref{eq:quadest} we get for any $r \in \mathbb{N}$ that
$$ \lambda_r \Norm{\Pi_{> r}(\phi_i \phi_j) }_{L^2}^2 \lesssim (\lambda_n + 1) \max_{1 \leq i \leq n}{ \|\phi_i\|^2_{L^{\infty}}}.$$

Using Weyl's law, we have (independently of the domain or manifold $\Omega$ and with an implicit constant converging to 1 as $\min{(r,n)} \rightarrow \infty$)
\begin{equation*}
  \lambda_r \sim \lambda_n \Bigl( \frac{r}{n} \Bigr)^{2/d}. 
\end{equation*}
and therefore
$$
  \Norm{\Pi_{> r}(\phi_i \phi_j) }_{L^2} \lesssim \Bigl( \frac{n}{r} \Bigr)^{1/d} \max_{1 \leq i \leq n} \norm{\phi_i}_{L^{\infty}}.
$$
Setting 
$$ r \sim  \left( \frac{1}{\varepsilon}\max_{1 \leq i \leq n} \norm{\phi_i}_{L^{\infty}} \right)^{d} n $$
makes the projection smaller than $\lesssim \varepsilon$.
\end{proof}

\begin{proof}[Proof of Corollary 2]
  To pass from Theorem 1 to Corollary 2 is straightforward: we make
  use of the $L^{\infty} - L^2$ estimate of elliptic eigenfunctions
  due to H\"ormander \cite{hor} (see Sogge \cite{sog} for the
  extension to general elliptic operators): let $\phi$ be an $L^2$-normalized
  eigenfunction of $L$ on a compact smooth manifold  associated with eigenvalue $\lambda$, we
  have
\begin{equation*}
  \norm{\phi}_{L^\infty} \leq C \lambda^{\frac{d - 1}{4}}. 
\end{equation*}
Combined with the Weyl's law, we thus obtain 
\begin{equation}
  \max_{1 \leq i \leq n} \norm{\phi_i}_{L^\infty} \lesssim_{\Omega, a, V} \lambda_n^{\frac{d - 1}{4}} \sim_{\Omega, a} \left(n^{2/d}\right)^{\frac{d - 1}{4}} \sim_{} n^{\frac{d-1}{2d}}.
\end{equation}
Corollary 2 immediately follows.
\end{proof}

\begin{proof}[Proof of Theorem 3]
We introduce the eigenfunctions of $-\Delta$ as $(\psi_k)_{k \in \mathbb{N}}$ having eigenvalues $\mu_k$. These eigenvalues are fairly closely
related to $\lambda_k$ which follows immediately from the min-max characterization
$$ \lambda_k = \inf_{V_k \subset C^{\infty}_c(\Omega)} \max_{ f \in V_k} \frac{ \int_{\Omega}{  |\nabla f|^2}}{ \int_{\Omega}{f^2}} \qquad \mbox{and} \qquad \mu_k = \inf_{V_k \subset C^{\infty}_c(\Omega)} \max_{ f \in V_k} \frac{ \int_{\Omega}{a(x)  |\nabla f|^2 + V f^2}}{ \int_{\Omega}{f^2}},$$
where the infimum runs over all $k-$dimensional subspaces of functions in $C^{\infty}_c(\Omega)$. This implies
$$\left( \min_{x \in \Omega} a(x) \right) \lambda_k - \|V\|_{L^{\infty}} \leq \mu_k \leq \left( \max_{x \in \Omega} a(x) \right) \lambda_k + \|V\|_{L^{\infty}}.$$
This shows that in the spectral definition of the $\dot H^{-1}-$norm
$$ \| f\|_{\dot H^{-1}}^2 = \sum_{k=1}^{\infty}{ \frac{ |\left\langle f, \psi_k\right\rangle|^2 }{\mu_k}}$$
it is possible to obtain an equivalent norm by replacing the eigenvalues of the Laplacian $\mu_k$ with the eigenvalues
of the Schr\"odinger operator $\lambda_k$.
We will now study the size of the projection onto the orthogonal space of the first $r$ eigenfunctions of the Laplacian 
\begin{equation*}
  \begin{aligned}
    \Norm{\Pi_{> r}(\phi_i \phi_j) }_{\dot{H}^{-1}}^2 
    & = \sum_{k=r+1}^{\infty}  \frac{ |\average{\phi_i \phi_j, \psi_k}|^2}{\mu_k} \\
&\lesssim_{\Omega, a, V}  \sum_{k=r+1}^{\infty}  \frac{ |\average{\phi_i \phi_j, \psi_k}|^2}{\lambda_k}\\
&\lesssim \frac{1}{\lambda_r} \sum_{k=r+1}^{\infty}   |\average{\phi_i \phi_j, \psi_k}|^2 = \frac{1}{\lambda_r} \| \Pi_{>r}(\phi_i \phi_j )\|^2_{L^2},
  \end{aligned}
\end{equation*}
where the projection $\Pi_{>r}$ is now onto the orthogonal complement
of the first $r$ Laplacian eigenfunctions.  We again use
$$       \norm{ \nabla (\phi_i \phi_j) }_{L^2} \leq \norm{\nabla \phi_i}_{L^2} \norm{\phi_j}_{L^{\infty}} + \norm{\nabla \phi_j}_{L^2} \norm{\phi_i}_{L^{\infty}}$$
and
$$ \|\nabla \phi_i\|_{L^2}^2 \lesssim_{a} \int_{\Omega}{ -\mbox{div}(a(x) \nabla \phi_i) \phi_i} \lesssim_{a} \lambda_n + \|V\|_{L^{\infty}}.$$
Then, however, the equivalence of the eigenvalues implies
\begin{align*}
  \lambda_r \| \Pi_{>r}(\phi_i \phi_j )\|^2_{L^2} &\leq \sum_{k=1}^{\infty}{ \lambda_k   |\average{\phi_i \phi_j, \psi_k}|^2 } \lesssim_{\Omega, a, V} \sum_{k=1}^{\infty}{ \mu_k   |\average{\phi_i \phi_j, \psi_k}|^2 } 
\end{align*}
and integration by part shows
\begin{align*}
\sum_{k=1}^{\infty}{ \mu_k   |\average{\phi_i \phi_j, \psi_k}|^2 } &= \left\langle \sum_{k=1}^{\infty}{ \mu_k   \average{\phi_i \phi_j, \psi_k} \psi_k} ,  \sum_{k=1}^{\infty}{   \average{\phi_i \phi_j, \psi_k} \psi_k} \right\rangle\\
&= \int_{-\Omega}{ - \Delta(\phi_i \phi_j) \phi_i \phi_j dx}\\
&= \norm{ \nabla (\phi_i \phi_j) }_{L^2}^2 \lesssim_{a,V}   
 \lambda_n  \max_{1 \leq i \leq n}{\| \phi_i \|_{L^{\infty}}^2}.
\end{align*}
Using Weyl's law, the above inequalities yield 
\begin{equation*}
  \begin{aligned}
  \Norm{\Pi_{> r}(\phi_i \phi_j) }_{\dot{H}^{-1}}  & \lesssim \frac{\sqrt{\lambda_n}}{\lambda_r} \max_{1 \leq i \leq n}{\| \phi_i \|_{L^{\infty}}} \\
&\lesssim n^{1/d} r^{-2/d}  \max_{1 \leq i \leq n}{\| \phi_i \|_{L^{\infty}}}.
  \end{aligned}
\end{equation*}
This is smaller than $\varepsilon$ for
$$ r \gtrsim \left( \frac{1}{\varepsilon}  \max_{1 \leq i \leq n}{\| \phi_i \|_{L^{\infty}}} \right)^{\frac{d}{2}} \sqrt{n}.$$
The unconditional result Corollary 4 follows from another application of 
$$ \max_{1 \leq i \leq n} \norm{\phi_i}_{L^\infty} \lesssim n^{\frac{d-1}{2d}}.$$
\end{proof}

\end{document}